\documentclass{gtart}

\def\ifplaintex{\expandafter\ifx\csname documentclass\endcsname\relax}

\expandafter\ifx\csname epsfbox\endcsname\relax\input epsf\fi

\ifplaintex 
\else
\fi

\def\gt{{\mathsurround=0pt\it $\cal G\mskip-2mu$eometry \&\ 
$\cal T\!\!$opology}}        

\def\gtp{{\mathsurround=0pt\it $\cal G\mskip-2mu$eometry \&\ 
$\cal T\!\!$opology $\cal P\!$ublications}}  

\def\lognumber#1{\def\thelognumber{#1}}
\def\volumenumber#1{\def\thevolumenumber{#1}}
\def\papernumber#1{\def\thepapernumber{#1}}
\def\volumeyear#1{\def\thevolumeyear{#1}}

\def\pagenumbers#1#2{\def\startpage{#1}\def\finishpage{#2}}
\def\published#1{\def\publishdate{#1}}
\def\proposed#1{\def\theproposer{#1}}
\def\seconded#1{\def\theseconders{#1}}
\def\received#1{\def\receiveddate{#1}}
\def\revised#1{\def\reviseddate{#1}}
\def\accepted#1{\def\accepteddate{#1}}

\def\asciiaddress#1{\def\theasciiaddress{#1}}
\def\asciiemail#1{\def\theasciiemail{#1}}
\long\def\asciiabstract#1{\long\def\theasciiabstract{#1}}
\def\asciikeywords#1{\def\theasciikeywords{#1}}

\let\\\par\let\thevolumenumber\relax\let\thepapernumber\relax
\let\thevolumeyear\relax\let\thesamplenumber\relax\let\startpage\relax
\let\finishpage\relax\let\publishdate\relax\let\receiveddate\relax
\let\reviseddate\relax\let\accepteddate\relax\let\theasciititle\relax
\let\theasciiauthors\relax\let\theasciiaddress\relax
\let\theasciiabstract\relax\let\theasciikeywords\relax
\let\theasciiemail\relax\let\theshortauthors\relax\let\theshorttitle\relax

\long\def\maketitlep{   

\count0=\startpage

\gt\hfill      
\hbox to 77pt{\vbox to 0pt{\vglue -15pt\epsfbox{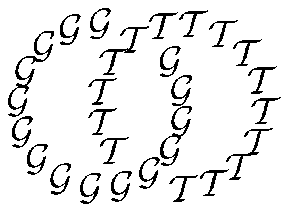}\vss}\hss}
\break
{\small\ifx\thesamplenumber\relax 
Volume \else Sample
\fi\thevolumenumber\ (\thevolumeyear)
\startpage--\finishpage\nl
Published: \publishdate}
\vglue 0.5truein plus 0.4fil minus 0.1truein

{\parskip=0pt\leftskip 0pt plus 1fil\def\\{\par\smallskip}{\ifplaintex\large
\else\Large\fi\bf\thetitle}\par\medskip}   

\vglue 0pt plus 0.1fil 

{\parskip=0pt\leftskip 0pt plus 1fil\def\\{\par}{\sc\theauthors}
\par\medskip}

\vglue 0pt plus 0.1fil 

{\small\parskip=0pt\let\newline\\
{\leftskip 0pt plus 1fil\def\\{\par}{\sl\theaddress}\par}
\expandafter\ifx\theemail\relax    
\relax\else\vglue 5pt plus 0.02fil minus 2pt\def\\{\stdspace{\rm 
and}\stdspace} 
\cl{Email:\stdspace\tt\theemail}\fi
\ifx\theurl\relax                  
\relax\else\vglue 5pt plus 0.02fil minus 2pt\def\\{\stdspace{\rm 
and}\stdspace}
\cl{URL:\stdspace\tt\theurl}\fi\par}

\vglue 7pt plus 0.3fil minus 3pt

{\bf Abstract}
\vglue 5pt plus 0.1fil minus 2pt

\theabstract

\vglue 7pt plus 0.3fil minus 3pt

{\bf AMS Classification numbers}\quad Primary:\quad \theprimaryclass

Secondary:\quad \thesecondaryclass

\vglue 5pt plus 0.3fil minus 2pt

{\bf Keywords:}\quad \thekeywords

\vglue 10pt plus 0.5fil minus 5pt

{\small  Proposed: \theproposer\hfill Received: \receiveddate\nl
Seconded: \theseconders\hfill 
\ifx\reviseddate\relax                         
Accepted: \accepteddate                        
\else
Revised: \reviseddate                          
\fi}
\eject
}       

\let\maketitlepage\maketitlep
\let\maketitle\maketitlepage

\font\phead=cmsl9 scaled 950
\font=cmsl9 scaled 1050
\font\pnum=cmbx10 scaled 913
\font\lnum=cmbx10 
\font\pfoot=cmsl9 scaled 950
\font=cmsl9 scaled 1050
\ifplaintex
\headline{\vbox to 0pt{\vskip -4.5mm\line{\small\phead\ifnum
\count0=\startpage ISSN 1364-0380 (on line)
1465-3060 (printed) \hfill {\pnum\folio}\else\ifodd\count0\def\\{ }%
\ifx\theshorttitle\relax\thetitle\else\theshorttitle\fi\hfill{\pnum\folio}
\else\def\\{ and }{\pnum\folio}\hfill\ifx\theshortauthors\relax\theauthors
\else\theshortauthors\fi\fi\fi}\vss}}
\footline{\vbox to 0pt{\vglue 0mm\line{\small\pfoot\ifnum\count0=\startpage
\copyright\ \gtp\hfill\else
\gt, Volume \thevolumenumber\ (\thevolumeyear)\hfill\fi}\vss
}}
\else
\makeatletter
\def\@oddhead{{\small\lhead\ifnum\count0=\startpage ISSN 1364-0380 (on line)
1465-3060 (printed) \hfill {\lnum\number\count0}\else\ifodd\count0
\def\\{ }\ifx\theshorttitle\relax \thetitle \else\theshorttitle\fi\hfill
{\lnum\number\count0}\else\def\\{ and }{\lnum\number\count0}
\hfill\ifx\theshortauthors\relax 
\theauthors\else\theshortauthors\fi\fi\fi}}\def\@evenhead{@oddhead}
\def\@oddfoot{\small\lfoot\ifnum\count0=\startpage\copyright\ \gtp\hfill\else
\gt, Volume \thevolumenumber\ (\thevolumeyear)\hfill\fi}
\def\@evenfoot{@oddfoot}
\makeatother
\fi

\newwrite\gtoutfile
\long\gdef\makeheadfile{  
{\def\\{, }\def\s{ }
\immediate\openout\gtoutfile head.xxx
\immediate\write\gtoutfile{To: math@arxiv.org}
\immediate\write\gtoutfile{Subject: put OR rep NNNNN:pppp}
\immediate\write\gtoutfile{--text follows this line--}
\immediate\write\gtoutfile{Proxy-for: \ifx\theasciiauthors\relax
\theauthors\else\theasciiauthors\fi\s<\ifx\theasciiemail\relax\theemail\else\theasciiemail\fi>}
\immediate\write\gtoutfile{\noexpand\\}
\immediate\write\gtoutfile{Authors: \ifx\theasciiauthors\relax
\theauthors\else\theasciiauthors\fi}
{\def\\{ }\immediate\write\gtoutfile{Title: \ifx\theasciititle\relax
\thetitle\else\theasciititle\fi}}
\immediate\write\gtoutfile{Subj-class: GT or GR or SG or ...}
\immediate\write\gtoutfile{MSC-class: \theprimaryclass\ifx\thesecondaryclass\relax\else, \thesecondaryclass\fi}
\immediate\write\gtoutfile{Journal-ref: Geom. Topol. \thevolumenumber\s
(\thevolumeyear) \startpage-\finishpage}
\immediate\write\gtoutfile{Comments: Published in Geometry and Topology at}
\immediate\write\gtoutfile{    http://www.maths.warwick.ac.uk/gt/GTVol\thevolumenumber/paper\thepapernumber.abs.html}
\immediate\write\gtoutfile{\noexpand\\}
\immediate\write\gtoutfile{}
\ifx\theasciiabstract\relax
\immediate\write\gtoutfile{\theabstract}\else
\immediate\write\gtoutfile{\theasciiabstract}\fi
\immediate\write\gtoutfile{}
\immediate\write\gtoutfile{\noexpand\\}
\immediate\write\gtoutfile{}
\immediate\closeout\gtoutfile}}  

\def\maketitlepage{\maketitlep\makeheadfile}
\let\maketitle\maketitlepage

\lognumber{184}
\volumenumber{6}\papernumber{3}\volumeyear{2002}
\pagenumbers{59}{67}
\received{24 May 2001}
\accepted{26 February 2002}
\revised{7 February 2002}
\proposed{Dieter Kotschick}
\seconded{Walter Neumann, Gang Tian}
\published{27 February 2002}

\usepackage{diagrams}
\diagramstyle{height=.25in}

\usepackage{amsmath,amsfonts}

\newtheorem{thm}{Theorem}[section]

\newtheorem{cor}[thm]{Corollary}

\newtheorem{problem}{Problem}
\newtheorem{rem1}[thm]{Remark}

\newcommand{\rnums} {{\mathbf R}}		
\newcommand{\znums} {{\mathbf Z}}		

\newcommand{\Ker}{\operatorname{Ker}}

\newcommand{\til}[1]{\Tilde{#1}}

\newcommand{\M}{{\mathcal{M}}}

\newcommand{\Norm}{\operatorname{Nm}}
\newcommand{\Picnaught}{\operatorname{Pic}^{0}}

\renewcommand{\tilde}[1]{\widetilde{#1}}
\renewcommand{\Tilde}[1]{\widetilde{#1}}

\begin{document}

\title{Surface bundles over surfaces of small genus}
\author{Jim Bryan\\Ron Donagi}

\address{
Department of Mathematics, University of British Columbia\\
121-1984 Mathematics Road, Vancouver BC 
\\Canada V6T 1Z2\\
}
\secondaddress{
Department of Mathematics, University of Pennsylvania\\
209 S 33rd Street, Philadelphia, PA 19104-6395, USA
}

\asciiaddress{
Department of Mathematics, University of British Columbia\\
121-1984 Mathematics Road, Vancouver BC 
\\Canada V6T 1Z2\\and\\Department of Mathematics, University of Pennsylvania\\
209 S 33rd Street, Philadelphia, PA 19104-6395, USA
}

\asciiemail{jbryan@math.ubc.ca, donagi@math.upenn.edu}
\email{jbryan@math.ubc.ca}\secondemail{donagi@math.upenn.edu}

\begin{abstract}
We construct examples of non-isotrivial algebraic families of smooth
complex projective curves over a curve of genus 2. This solves a problem
from Kirby's list of problems in low-dimensional topology. Namely, we show
that 2 is the smallest possible base genus that can occur in a 4--manifold
of non-zero signature which is an oriented fiber bundle over a Riemann
surface. A refined version of the problem asks for the minimal base genus
for fixed signature and fiber genus. Our constructions also provide new
(asymptotic) upper bounds for these numbers.
\end{abstract}

\asciiabstract{
We construct examples of non-isotrivial algebraic families of smooth
complex projective curves over a curve of genus 2. This solves a problem
from Kirby's list of problems in low-dimensional topology. Namely, we show
that 2 is the smallest possible base genus that can occur in a 4-manifold
of non-zero signature which is an oriented fiber bundle over a Riemann
surface. A refined version of the problem asks for the minimal base genus
for fixed signature and fiber genus. Our constructions also provide new
(asymptotic) upper bounds for these numbers.}

\asciikeywords{Surface bundles, 4-manifolds, algebraic surface}



\maketitle



\section{Introduction}

By a \emph{surface bundle over a surface} we will mean an oriented fiber
bundle whose fibers are compact, oriented 2--manifolds and whose base is a
compact, oriented 2--manifold.  In this paper, we solve the following problem,
posed by Geoff Mess, from Kirby's problem list in low-dimensional topology:
\begin{problem}[Mess, \cite{Kirby-problems} Problem 2.18A] What is the
smallest number $b$ for which there exists a surface bundle over a surface
with base genus $b$ and non-zero signature?
\end{problem}

The first examples of surface bundles over surfaces with non-zero signature
were constructed independently by Atiyah \cite{Atiyah-example} and Kodaira
\cite{Kodaira-example} (which were then generalized by Hirzebruch in
\cite{hirzebruch-example}); these examples had base genus 129. In his
remarks following the statement of the problem, Mess alludes to having a
construction with base genus 42; later examples with base genus 9 were
constructed in \cite{EKKOS}. Subsequently, it was noticed by several people
(eg \cite{BDS,LeBrun}) that the original examples of Atiyah, Kodaira, and
Hirzebruch have two different fibrations, one of which is over a surface of
genus 3.

Since the signature of a 4--manifold which fibers over a sphere or torus
must vanish, the smallest possible base genus is two. We prove that this
does indeed occur as a special case of our main construction.

\begin{thm}\label{thm: main construction}
For any integers $g,n\geq 2$, there exists a connected algebraic surface
$X_{g,n}$ of signature $\sigma (X_{g,n})=\frac{4}{3}g (g-1)
(n^{2}-1)n^{2g-3}$ that admits two smooth fibrations $\pi _{1}\co X_{g,n}\to
C$ and $\pi _{2}\co X_{g,n}\to \tilde{D}$ with base and fiber genus
$(b_{i},f_{i})$ equal to
\begin{align*}
(b_{1},f_{1})&= (g,g(gn-1)n^{2g-2} +1)\text{ and}\\
(b_{2},f_{2})&= (g(g-1) n^{2g-2}+1,gn)
\end{align*}
respectively.
\end{thm}

In particular, for $n=g=2$ the manifold $X_{2,2}$ from Theorem~\ref{thm:
main construction} gives us:

\begin{cor}
There exists a 4--manifold of signature 16 that fibers over a
surface of genus 2 with fiber genus 25.
\end{cor}

Any surface bundle $X\to B$ with fiber genus $f$ is determined up to
isomorphism by the homotopy class of its classifying map $\phi \co B\to \M
_{f}$, where $\M _{f} $ is the moduli space of non-singular genus $f$
curves, regarded as a complex orbifold, and $\phi $ is an orbi-map (and the
homotopy class is formed using homotopies in the orbifold category).

From the index theorem for families (see \cite{Atiyah-example} or
\cite{Ivan-Smith}), the signature of $X$ is determined by the evaluation of
 the first Chern class of the Hodge bundle $\mathbb{E}\to \M _{f}$ on $B$:
\[
\sigma (X)=4\int _{B}\phi ^{*} (c_{1} (\mathbb{E})).
\]
Since for $f\geq 3$, $\det (\mathbb{E})$ is ample on $\M _{f}$
(eg \cite{Harris-Morrison-book}), $\phi ^{*} (c_{1} (\mathbb{E}))$ will
evaluate non-trivially on $B$ for any non-constant holomorphic orbi-map
$\phi \co B\to \M _{f}$. Thus any holomorphic family $X\to B$ that is not
isotrivial will have non-zero signature.

For $f\geq 3$, the non-torsion part of $H_{2} (\M _{f};\znums )$ is of rank
one and is generated by the dual of $c_{1} (\mathbb{E})$ and so one can
refine the original problem as the problem of determining the minimal genus
for representatives of elements of $H_{2} (\M _{f};\znums )$ mod torsion
(c.f. \cite{Kirby-problems} 2.18B and \cite{EKKOS}). That is, one can try
to find the numbers:
\[
b_{f} (m)=\min\{b:\text{ $\exists $ a genus $f$ bundle $X\to B$ with $g
(B)=b$ and $\sigma (X)=4m$.} \}
\]
Kotschick has determined lower bounds on $b_{f} (m)$ using Seiberg--Witten
theory \cite{Kotschick}, and the constructions of \cite{Endo} and later
\cite{EKKOS} give systematic upper bounds for $b_{f} (m)$.  Given a bundle
$X\to B$, one obtains a sequence of bundles by pulling back by covers of
the base. The base genus and signature grow linearly in this sequence, so
it is natural to consider the minimal genus asymptotically. Define
\[
G_{f}=\lim_{m\to \infty }\frac{b_{f} (m)}{m}.
\]
It is easy to see that this limit exists and is finite (see
\cite{Kirby-problems} 2.18B). Upper bounds for $G_{f}$ are given by Endo,
et al in \cite{EKKOS}; our constructions substantially improve their upper
bounds for the case when $f$ is composite:
\begin{cor}
Let $G_{f}$ be defined as above and suppose that $f=ng$ with $n,g\geq 2$. Then
\[
G_{f}\leq \frac{3n}{n^{2}-1}.
\]
\end{cor}

\proof Start with the bundle $X_{g,n}\to \til{D}$ from the theorem
and construct a sequence of bundles $X_{g,n}^{m}\to \til{D}^{m}$ obtained
by pulling back by unramified, degree $m$ covers of the base
$\til{D}^{m}\to \til{D}$. The signature and base genus of these examples
are easily computed:
\begin{align*}
\sigma (X_{g,n}^{m})&=m\sigma (X_{g,n})\\
g (\til{D}^{m})-1&=m (g (\til{D})-1)
\end{align*}
and so 
$$
G_{f}\leq \lim_{m\to \infty }\frac{mg (g-1)n^{2g-2}+1}{\tfrac{m}{3}g (g-1) (n^{2}-1)n^{2g-3}}=\frac{3n}{n^{2}-1}.\eqno{\qed}
$$
\medskip

For example, if 
$f$ is even, then we have 
\[
G_{f}\leq \frac{6f}{f^{2}-4}<\frac{6}{f-2}
\]
which improves the bound of $\frac{16}{f-2}$ found in \cite{EKKOS}. Note
that Kotschick's lower bound is $\frac{2}{f-1}$.

Our constructions are similar to Hirzebruch, Atiyah, and Kodaira's in that
they are also branched covers of a product of Riemann surfaces. We have
refined and extended their approach and we also employ some ideas that go
back to a construction of Gonzalez-Diez and Harvey \cite{Diez-Harvey}. We
would like to thank Dieter Kotschick for helpful comments and suggestions.

The first author is supported by an Alfred P Sloan Research
Fellowship and NSF grant DMS-0072492 and the second author is supported by
NSF grant DMS-9802456.

\section{The main construction}\label{sec: main construction}

We will construct $X_{g,n}$ as a degree $n$, cyclic branched cover of a
certain product of curves, $\til{D}\times C$. This cover will be branched
along two disjoint curves $\Gamma _{1}$ and $\Gamma _{2}$ where the $\Gamma
_{i}$'s are the graphs of unramified maps $\til{f}_{i}\co \til{D}\to C$. We
begin by first constructing intermediate covers $f_{i}\co D\to C$.

We construct $D$ and $C$ as follows. Fix an elliptic curve $E$ with origin
$o\in E$ and fix a 2--torsion point $\tau \in E$. Let $\pi \co C\to E$ be a
$g$--fold cyclic cover of $E$ branched at $o$ and $\tau $.
Note that the genus of $C$ is $g$. Let $x\mapsto x+\tau $ denote
translation by $\tau $. We define $D'\subset C\times C$ to be the locus of
points $(p_{1},p_{2})$ such that $\pi (p_{1})=\pi (p_{2})+\tau $. $D'$ is
clearly disjoint from the diagonal and $D'$ has two maps $f'_{i}\co D'\to C$
induced by the projections. Consider the preimage of a point $p_{1}\in C$
under the map $f_{1}'$. It is all pairs of the form $(p_{1},\pi ^{-1} (\pi
(p_{1})+\tau ))$ and so $f_{i}'$ is of degree $g$ and is unramified away
from the two points $(\pi ^{-1} (o),\pi ^{-1} (\tau ))$ and $(\pi ^{-1}
(\tau ),\pi ^{-1} (o))$. We will show that these points are ordinary
$g$--fold singularities of $D'$ and so then letting $D\to D'$ be the
normalization, we will obtain the unramified, degree $g$ covers $f_{i}\co D\to
C$ by the composition of $f_{i}'$ with the normalization.

To see that $(\pi ^{-1} (o),\pi ^{-1} (\tau ))\in D'$ is an ordinary
$g$--fold singular point, consider local coordinates $u$ and $v$ on $E$
about $o$ and $\tau $ such that $u$ is identified to $v$ by translation by
$\tau $. Choose local coordinates $z$ and $w$ on $C$ so that $\pi $ is
locally given by $u=z^{g} $ and $v=w^{g}$. Then $z^{g}=w^{g}$ are the local
equations for $D'$ in $C\times C$ at the points $(\pi ^{-1} (o),\pi ^{-1}
(\tau ))$ and $(\pi ^{-1} (\tau ),\pi ^{-1} (o))$ which are thus ordinary
$g$--fold singularities.

Note that since $D'$ is disjoint from the diagonal, the covers $f_{i}\co D\to
C$ have the property that $f_{1} (p)\neq f_{2} (p)$ for all $p\in D$. It is
not immediately clear from the construction that $D$ is connected; we will
postpone the discussion of this issue until the end of the section.

We next construct the unramified cover $\til{D}\to D$.  Let $\Norm
\co \Picnaught (C)\to \Picnaught (E) $ be the norm map induced by $\pi $, that
is, given a degree zero divisor $\sum m_{i}p_{i}$ on $C$, $\Norm (\sum
m_{i}p_{i})$ is defined by $\sum m_{i}\pi (p_{i})$. Note that by
construction, 
$$\Norm (\mathcal{O} (p_{1}-p_{2}))=\mathcal{O} (\tau -o)\quad\text{for}\quad
(p_{1},p_{2})\in D'\subset C\times C.$$ 
We choose an $n$th root of
$\mathcal{O} (\tau -o)$ which we denote by $R$.

We define an unramified cover $\tilde{D}\to D$ of degree $n^{2g-2}$ as
follows.
\[
\tilde{D}=\left\{(L,(p_{1},p_{2}))\in  \Picnaught (C)\times D:\quad
L^{\otimes n}\cong \mathcal{O} (p_{1}-p_{2}),\quad \Norm (L)=R
\right\}.
\]
The natural projection $\tilde{D}\to D$ is unramified and has degree
$n^{2g-2}$ since the fibers are torsors on the $n$--torsion points in $\Ker
(\Norm )$ (which is a connected Abelian variety of dimension $g-1$ by the
argument below). Let $\tilde{f}_{i}\co \tilde{D}\to C$ be the compositions
with $f_{i}$ and let $\Gamma _{i}\subset \tilde{D}\times C$ be the
corresponding graphs. Since $\tilde{f}_{1} (\tilde{p})\neq \tilde{f}_{2}
(\tilde{p})$ for all $\tilde{p}\in \tilde{D}$, the curves $\Gamma _{1}$ and
$\Gamma _{2}$ are disjoint.  We will discuss the connectedness of
$\tilde{D}$ at the end of the section.

To see that $\Ker (Nm)$ is connected, consider the following diagram with
exact rows:
\[
\begin{diagram}
0&\rTo &H_{1} (C,\znums )&\rTo&H_{1} (C,\rnums ) &\rTo &\Picnaught (C)&\rTo&0\\
 &     &\dTo_{a_{1}}     &    & \dTo_{a_{2}}     &     &\dTo_{\Norm }     &    & \\
0&\rTo &H_{1} (E,\znums )&\rTo&H_{1} (E,\rnums ) &\rTo &\Picnaught (E)&\rTo&0
\end{diagram}
\]

Since $\Ker (a_{2})$ is connected, $\Ker (\Norm )$ is connected if $\Ker
(a_{2})\to \Ker (\Norm )$ is surjective. By a diagram chase, $\Ker
(a_{2})\to \Ker (\Norm )$ is surjective if $a_{1}$ is surjective. But
$a_{1}$, which is $\pi _{*}$, is indeed surjective because $\pi $ does not
factor through any unramified cover (the factored map would have to have
only one ramification point which is impossible).

We want to construct $X_{g,n}\to \tilde{D}\times C$ as a cyclic branched
cover of degree $n$, ramified over $\Gamma _{1}-\Gamma _{2}$. To do this we
need to construct a line bundle $\mathcal{L}\to \tilde{D}\times C$ so that
$\mathcal{L}^{\otimes n}\cong \mathcal{O} (\Gamma _{1}-\Gamma _{2})$. Once
we have $\mathcal{L}$, we will define
\[
X_{g,n}=\{(v_{1}:v_{2})\in \mathbf{P} (\mathcal{L}\oplus \mathcal{O}):\quad
(v^{n}_{1}:v^{n}_{2})= (s_{1}:s_{2})\}
\]
where $s_{i}$ is a section of $\mathcal{O} (\Gamma _{i})$ that vanishes
along $\Gamma _{i}$ so that $(s_{1}:s_{2})$ is in $  \mathbf{P}
(\mathcal{O} (\Gamma _{1})\oplus \mathcal{O} (\Gamma _{2})) $ which is the
same as $\mathbf{P} (\mathcal{O} (\Gamma _{1}-\Gamma _{2})\oplus
\mathcal{O} ) $.

To find $\mathcal{L}$, we use the Poincar\'e bundle $\mathcal{P}\to \Picnaught
(C)\times C$ which is a tautological bundle in the sense that
$\mathcal{P}|_{\{L \}\times C}\cong L$. $\mathcal{P}$ is uniquely
determined by choosing a point $p_{0}\in C$ and specifying that
$\mathcal{P}$ restricted to $\Picnaught (C)\times \{p_{0} \}$ is
trivial. We use the same letter $\mathcal{P}$ to denote the pullback of
$\mathcal{P}$ by the composition of the inclusion and projection:
\[
\tilde{D}\times C\to \Picnaught (C)\times D\times C\to \Picnaught (C)\times C.
\]
Let $M\in \Picnaught (\tilde{D})$ be an $n$-th root of $\mathcal{O}
(\til{f}_{1}^{-1} (p_{0})-\til{f}_{2}^{-1} (p_{0}))$ and define
$\mathcal{L}$ to be $\mathcal{P}\otimes \pi ^{*}_{\Tilde{D}}M$. We need to
show that $\mathcal{L}^{\otimes n}\cong \mathcal{O} (\Gamma _{1}-\Gamma
_{2})$ or equivalently, $(\mathcal{L}^{\vee })^{\otimes n}\otimes
\mathcal{O} (\Gamma _{1}-\Gamma _{2})\cong \mathcal{O}$. Let $x=
(L,p_{1},p_{2})$ be any point of $\Tilde{D}$. By construction, we have
\begin{align*}
\mathcal{L}^{\otimes n}|_{\{x \}\times C}&\cong \mathcal{P}^{\otimes n}|_{\{x \}\times C}\\
&\cong L^{\otimes n}\\
&\cong \mathcal{O} (p_{1}-p_{2})\\
&\cong \mathcal{O} (\Gamma _{1}-\Gamma _{2})|_{\{x \}\times C} ,
\end{align*}
therefore, $(\mathcal{L}^{\vee })^{\otimes n}\otimes \mathcal{O} (\Gamma
_{1}-\Gamma _{2})$ is trivial on every slice $\{x \}\times C$ and so it
must be the pullback of a line bundle on $\Tilde{D}$. But
\begin{align*}
\mathcal{L}^{\otimes n}|_{\Tilde{D}\times p_{0}}&\cong \mathcal{P}^{\otimes n}|_{\Tilde{D}\times p_{0}}\otimes M^{\otimes n}\\
&\cong \mathcal{O} (\til{f}_{1}^{-1} (p_{0})-\til{f}_{2}^{-1} (p_{0}))\\
&\cong \mathcal{O} (\Gamma _{1}-\Gamma _{2})|_{\Tilde{D}\times p_{0}}
\end{align*}
and so $(\mathcal{L}^{\vee })^{\otimes n}\otimes \mathcal{O} (\Gamma
_{1}-\Gamma _{2})$ is indeed the trivial bundle.  The line bundle
$\mathcal{L}$ then gives us the $n$--fold cyclic branched cover $X_{g,n}\to
\Tilde{D}\times C$ by the construction described above.

The fiber of the projection $X_{g,n}\to \til{D}$ over a point $x=
(L,p_{1},p_{2})\in \til{D}$ is the $n$--fold cyclic branched cover of $C$
branched at $p_{1}-p_{2}$ determined by $L$. By the Riemann--Hurwitz
formula, this curve has genus $gn$. On the other hand, the fiber of
$X_{g,n}\to C$ over a point $p\in C$ is an $n$--fold cyclic cover of
$\til{D}$ branched over $\til{f}_{1}^{-1} (p)-\til{f}_{2}^{-1} (p)$ which
consists of $2gn^{2g-2}$ (distinct) points. Noting that $g (\til{D})=
g(g-1) n^{2g-2}+1$, one easily computes the fiber genus to be
$g(gn-1)n^{2g-2}+1$.

To determine the signature of $X_{g,n}$ we use a formula for the signature of
a cyclic branched cover due to Hirzebruch \cite{hirzebruch-example}:
\begin{equation}\label{eqn: hirzebruch's sign formula}
\sigma (X_{g,n})=\sigma (\til{D}\times C)-\frac{n^{2}-1}{3n} (\Gamma
_{1}-\Gamma _{2})^{2}.
\end{equation}
The signature of $\til{D}\times C$ is zero, and since $\Gamma _{1}$ and
$\Gamma _{2} $ are disjoint, we just need to compute $\Gamma
_{1}^{2}=\Gamma _{2}^{2}$. By the adjunction formula, we have
\begin{align*}
\Gamma _{i}^{2}&=2g (\til{D})-2 -K_{\til{D}\times C}\cdot \Gamma _{i}\\
&=2g (\til{D})-2 -\left(2g (\til{D})-2 + \operatorname{deg} (\til{f}_{i}) (2g (C)-2)\right)\\
&= -\operatorname{deg} (\til{f}_{i}) (2g (C)-2)\\
&=-2g (g-1)n^{2g-2}
\end{align*}
and so
\[
\sigma  (X_{g,n})=\frac{4}{3} g(g-1)  (n^{2}-1)n^{2g-3}.
\]
We have not yet proved that $X_{g,n}$ is connected since it is not clear
from their constructions whether $D$ and $\tilde{D}$ are connected or
not. If $D$ or $\til{D}$ were not connected, it would actually improve our
construction in the sense that the connected components of $X_{g,n}$ would
still fiber as surface bundles in two different ways but would have a
smaller base or fiber genus (depending on which fibration is
considered). In fact, for certain choices of $C$, one can show that $D$ is
disconnected when $g$ is a composite number with an odd factor. However, we
do not explore these possibilities but instead, to complete the proof of
Theorem~\ref{thm: main construction} as stated, we show that one can always
take $X_{g,n} $ to be connected.

To this end, suppose that $\til{D}$ is disconnected with $N$
components. Since $\til{D}\to D$ and $D\to C$ are normal coverings, $N$
must divide $gn^{2g-2}$, the degree of $\til{f}_{i}\co \til{D}\to C$. Fix a
connected component $\til{D}'$ of $\til{D}$ and let $X_{g,n}'$ be the
corresponding component of $X_{g,n}$. Note that $X_{g,n}'\to \til{D}'\times
C$ is the cyclic branched cover determined by
$\mathcal{L}':=\mathcal{L}|_{\til{D}'\times C}$. Note that the degree of
$\til{D}'\to C$ is $N^{-1}gn^{2g-2}$. Now consider any connected,
unramified, degree $N$ cover $p\co D''\to \til{D}'$ and let $f_{i}''\co D''\to C$
be the composition of $p$ with $\til{f}_{i}|_{\til{D}'}$ noting that the
degree of $f_{i}''$ is $gn^{2g-2}$. Let $\Gamma ''_{i}\subset D''\times C$
be the graph of $f_{i}''$ and observe that $p^{*} (\mathcal{L}')^{\otimes
n}\cong \mathcal{O} (\Gamma _{1}''-\Gamma _{2}'')$ so that $p^{*}
(\mathcal{L}')$ defines an $n$--fold cyclic branched cover $X''_{g,n}\to
D''\times C$ ramified along $\Gamma _{1}''-\Gamma _{2}''$.

The computation of the signature of $X_{g,n}''$ and the computation of the
base and fiber genera of the fibrations $X_{g,n}''\to D''$ and
$X_{g,n}''\to C$ then proceed identically with the corresponding
computations for $X_{g,n}$ done previously (where we were implicitly
assuming that $\til{D}$ was connected). Indeed, those computations only
depended upon the degree of $\til{f}_{i}$ which is the same as the degree
of $f_{i}''$. Therefore, whenever $\til{D}$ is not connected, we replace
$\til{D}$ with $D''$ and we replace $X_{g,n}$ with the connected surface
$X''_{g,n}$ thus completing the proof of Theorem~\ref{thm: main
construction}.

\subsection{A simple construction of a base genus 2 surface bundle}
The surfaces $X_{g,n}$ were constructed to be economical with both the
fiber genus and the base genus. A simple construction of a base genus 2
surface bundle (but with larger fiber genus) can be obtained as
follows. Let $C$ be a genus 2 curve with a fixed point free automorphism
$\sigma \co C\to C$ (eg, let $C$ be the smooth projective model of
$y^{2}=x^{6}-1$ which has a fixed point free automorphism of order 6 given
by $(x,y)\mapsto (e^{2\pi i/6}x,-y)$). Let $\rho \co \tilde{C}\to C$ be the
unramified cover corresponding to the surjection $\pi _{1} (C)\to H_{1}
(C,\znums /2)$. Then the graphs $\Gamma _{\rho }$ and $\Gamma _{\sigma
\circ \rho }$ are disjoint in $\til{C}\times C$ and the class $[\Gamma
_{\rho }]+[\Gamma _{\sigma \circ \rho }]$ is divisible by 2 (by an argument
similar to the one in \cite{BDS} for example). Therefore, there exists a
double cover, $X\to \til{C}\times C$ branched along $\Gamma _{\rho }$ and
$\Gamma _{\sigma \circ \rho }$, so that the projections $X\to C$ and $X\to
\til{C}$ are smooth fibrations. One then easily computes that the bundle
$X\to C$ has base genus 2, fiber genus 49, and signature 32.


\begin{thebibliography}

\bibitem{Atiyah-example}
\textbf{M\,F Atiyah}, \emph{The signature of fibre-bundles}, from: ``Global
  Analysis (Papers in Honor of K. Kodaira)'', Univ. Tokyo Press, Tokyo (1969)
  73--84

\bibitem{BDS}
\textbf{Jim Bryan{\rm,} Ron Donagi{\rm,} Andras Stipsicz}, \emph{{Surface bundles: some
  interesting examples}}, Turkish J. Math. 25 (2001) 61--68, proceedings of the
  {$7^{th}$ G\"okova} Geometry and Topology conference

\bibitem{EKKOS} \textbf{H Endo{\rm,} M Korkmaz{\rm,} D Kotschick{\rm,}
B Ozbagci{\rm,} A Stipsicz}, \emph{{Commutators, Lefschetz fibrations
and the signatures of surface bundles}}, {\tt arXiv:math.GT/0103176}, to
appear in Topology

\bibitem{Endo}
\textbf{Hisaaki Endo}, \emph{A construction of surface bundles over surfaces
  with non-zero signature}, Osaka J. Math. 35 (1998) 915--930

\bibitem{Diez-Harvey}
\textbf{Gabino Gonz{\'a}lez-D{\'\i}ez, William~J Harvey}, \emph{On complete
  curves in moduli space. {I}, {I}{I}}, Math. Proc. Cambridge Philos. Soc. 110
  (1991) 461--466, 467--472

\bibitem{Harris-Morrison-book} \textbf{Joe Harris{\rm,} Ian Morrison},
\emph{Moduli of curves}, Springer--Verlag, New York (1998)

\bibitem{hirzebruch-example}
\textbf{F Hirzebruch}, \emph{The signature of ramified coverings}, from:
  ``Global Analysis (Papers in Honor of K. Kodaira)'', Univ. Tokyo Press, Tokyo
  (1969)  253--265

\bibitem{Kirby-problems}
\textbf{Rob Kirby}, \emph{Problems in low dimensional topology}, from:
  ``Proceedings of the 1993 Georgia International Topology Conference held at
  the University of Georgia, Athens, GA, August 2--13, 1993'', (William~H
  Kazez, editor), American Mathematical Society, Providence, RI (1997)

\bibitem{Kodaira-example}
\textbf{K Kodaira}, \emph{A certain type of irregular algebraic surfaces}, J.
  Analyse Math. 19 (1967) 207--215

\bibitem{Kotschick}
\textbf{D Kotschick}, \emph{Signatures, monopoles and mapping class groups},
  Math. Res. Lett. 5 (1998) 227--234

\bibitem{LeBrun}
\textbf{Claude LeBrun}, \emph{{Diffeomorphisms, symplectic forms, and Kodaira
  fibrations}}, Geom. Topol. 4 (2000) 451--456

\bibitem{Ivan-Smith}
\textbf{Ivan Smith}, \emph{Lefschetz fibrations and the {H}odge bundle}, Geom.
  Topol. 3 (1999) 211--233

\end{thebibliography}
\end{document}